\documentclass[12pt,twoside,draft]{amsart}
\usepackage{amssymb}
\usepackage{latexsym}
\usepackage{amsfonts}
\usepackage{amsmath}
\usepackage{amssymb}
\usepackage{portland}
\usepackage{times}

\newtheorem{theorem}{Theorem}[section]

\theoremstyle{definition}

\renewcommand{\leq}{\leqslant}
\renewcommand{\geq}{\geqslant}

\newcommand{\aut}[1]{{\sf Aut}(#1)}

\newcommand{\lcs}[2]{\gamma_{#2}\left(#1\right)}

\newcommand{\F}{\mathbb F}

\begin{document}

\title{A computer-based approach to the classification of nilpotent Lie algebras}
\author{Csaba Schneider}
\address{Informatics Laboratory\\
Computer and Automation Research Institute of the
Hungarian Academy of Sciences\\
1111 Budapest L\'agym\'anyosi u.~11.\\
Hungary}
\email{csaba.schneider@sztaki.hu\protect{\newline} {\it WWW:}
www.sztaki.hu/$\sim$schneider}

\begin{abstract}
We adapt the $p$-group generation algorithm to classify small-dimensional 
nilpotent Lie algebras over small fields. Using an implementation of this algorithm, we list 
the nilpotent Lie algebras of dimension at most~9 over $\F_2$ and those of dimension at most~7 over $\F_3$ and $\F_5$.  
\end{abstract}

\date{18 June 2004}

\maketitle

\section{Introduction}\label{intro}

The classification of $n$-dimensional 
nilpotent Lie algebras over a given field $\F$ is a very difficult 
problem even for relatively small $n$.
The aim of this article is to present a series of computer calculations
that led to the following theorem. 

\begin{theorem}\label{main}
The number of isomorphism types of
$6$-dimensional nilpotent Lie algebras is $36$ over $\F_2$, and
$34$ over $\F_3$ and $\F_5$.
The number of isomorphism types of
$7$-dimensional nilpotent Lie algebras is $202$ over $\F_2$, $199$ 
over $\F_3$, and $211$ over $\F_5$. 
The number of isomorphism types of
nilpotent Lie algebras with dimension $8$ and $9$ over $\F_2$ is $1831$ and $27073$, respectively.
\end{theorem}

The classifications in Theorem~\ref{main} were 
obtained via a series of computer calculations using a 
{\sf GAP~4}~\cite{gap} implementation of a nilpotent Lie algebra generation algorithm.
The ideas used in these calculations are the same as those used in the 
classification of finite $2$-groups with order at most $2^9$; see~\cite{p-gen, o'be}.
Let $\lcs Li$ denote the $i$-th term of the lower central series of a Lie algebra
$L$, so that $\lcs L1=L$, $\lcs L2=L'$, etc.
If $L$ is a finitely generated nilpotent Lie algebra 
with nilpotency class $c$ then $L$ is an immediate descendant of 
$L/\lcs L{c}$ (see Section~\ref{coversec} for definitions). 
Further, $L/\lcs Lc$ is an immediate descendant of $L/\lcs L{c-1}$. Continuing 
this way, we can see that every finitely generated nilpotent 
Lie algebra can be obtained after finitely many steps from a finite-dimensional
 abelian Lie 
algebra by computing immediate descendants. 
This suggests that a theoretical algorithm to generate all $n$-dimensional, 
nilpotent Lie algebras can be designed once we can efficiently compute
immediate descendants.
We will see that every immediate descendant of 
$L$ is a quotient of another nilpotent Lie algebra, which is referred 
to as the {\em
 cover}. It is shown in this paper that, for 
a finite-dimensional, nilpotent $\F_p$-Lie algebra $L$, 
it is possible to effectively compute
the  cover, and then to compute
a complete and irredundant list of the isomorphism types of 
the immediate descendants of $L$. Repeating
the immediate descendant calculation finitely many times, it is, 
in theory, possible to obtain a complete and irredundant list of 
all isomorphism types of the nilpotent Lie algebras 
with a given dimension over a finite field. In practice, this calculation quickly becomes unfeasible
as $n$ grows. Nevertheless, using this approach, it is possible to 
obtain classifications of Lie algebras that would otherwise be beyond
hope; see Theorem~\ref{main}.

The structure of this paper is as follows. 
In Section~\ref{coversec}
we develop the theory of a Lie algebra generation algorithm.
An application of the algorithm to prove Theorem~\ref{main}
will be presented in 
Section~\ref{class}. 
The final Section~\ref{implsec} will discuss an implementation 
of the algorithm.

\section{A nilpotent Lie algebra generation algorithm}\label{coversec}

Our nilpotent Lie algebra generation  algorithm is an 
adaptation of O'Brien's $p$-group generation algorithm, whose details 
can be found in~\cite{p-gen}. The Lie algebra 
generation algorithm is described without proofs in this section. Another 
variation on this theme is presented in~\cite{obnv-l} 
where the authors
classified  groups and nilpotent Lie rings of order $p^6$. Recently 
Michael Vaughan-Lee used the same approach to extend these results to $p^7$. 

Throughout this section $L$ is 
a finite-dimensional, 
nilpotent Lie algebra. Let $Z(L)$ denote the center of $L$. 
A nilpotent Lie algebra $K$ is said to be a {\em central extension} 
of $L$ if $K$ has an 
ideal $I$ such that $I\leq K'\cap Z(K)$ and $K/I\cong L$. In the terminology of~\cite{batten1,batten2}, $(K,I)$ is said to be a defining pair for $L$.
The algebra $K$ is said to be an
{\em immediate descendant} 
of $L$ if $L\cong K/\lcs Kc$ where $c$ is the nilpotency 
class of $K$.
Hence an immediate descendant is a special kind of central extension. 
The {\em cover} $L^*$ of a finitely generated nilpotent Lie algebra $L$ is defined as follows. Suppose 
that $\dim L/L'=d$. Then $L$ is a $d$-generator Lie algebra, and so the free 
Lie algebra $F_d$ with rank $d$ has an ideal $I$ such that $F_d/I\cong L$. 
The cover $L^*$ is defined as the Lie algebra $F_d/[I,F_d]$. The {\em multiplicator}
of $L^*$ is the ideal $I/[I,F_d]$. The cover 
$L^*$ is also a finite-dimensional 
nilpotent Lie algebra. Moreover, if $L$ has nilpotency class $c$ 
then the class of $L^*$ is at most $c+1$, and $\lcs {L^*}{c+1}$ is referred to 
as the {\em nucleus} of $L^*$. 

Suppose now without loss of generality that $L=F_d/I$ as in the previous
paragraph. Let $L^*$ be the cover of $L$ with 
multiplicator $M$ and nucleus $N$. 
Then $K$ is a central extension of $L$ if
and only if $K\cong L^*/J$ for some ideal $J\leq M$.
 Further, in this case, $K$ is an 
immediate descendant of $L$ if and only if $J\neq M$ and 
$J+N=M$. 
A proper subspace $J$ of $M$ with $J+N=M$ is said to
be {\em allowable}. Thus it is possible to obtain a complete list of immediate
descendants of $L$ by listing all quotients $L^*/J$ where $J$ runs through the 
allowable subspaces of the multiplicator $M$. Unfortunately, two different
allowable subspaces may lead to isomorphic Lie algebras. This problem can, 
however, be tackled using the automorphism group of $L$. If $\alpha$ is an 
automorphism, then $\alpha$ can be lifted to an automorphism $\alpha^*$ of the 
cover $L^*$ as follows. Let $\psi:L^*\rightarrow L$ denote the natural 
epimorphism with kernel $M$. Suppose that $b_1,\ldots,b_d$ is a minimal 
generating set for $L^*$; then $b_1\psi,\ldots,b_d\psi$ is a minimal generating set for $L$. Suppose that $y_1,\ldots,y_d\in L^*$ are 
chosen so that
$b_i\psi\alpha=y_i\psi$ for all $i\in\{1,\ldots,d\}$.  
Then the map $b_i\mapsto y_i$, for $i=1,\ldots,d$, 
can uniquely be extended to an automorphism of $L^*$. This automorphism 
is denoted $\alpha^*$, even though it is not uniquely determined by $\alpha$. 
On the other hand the restriction of $\alpha^*$ to $M=I/[I,F_d]$ depends only 
on $\alpha$. This defines a linear representation 
\begin{equation}\label{ro}
\varrho:\aut L\rightarrow
{\sf GL}(M)\quad\mbox{given by}\quad\alpha\mapsto\alpha^*|_M.
\end{equation}
Using a familiar argument, it is not hard to see that two 
allowable subspaces $J_1$ and $J_2$ give isomorphic Lie algebras $L^*/J_1$ 
and $L^*/J_2$ if and only if $J_1$ and $J_2$ are in the same orbit under the 
action $\aut L\varrho$. 

If $J$ is an allowable subspace of the multiplicator then the automorphism group
of $K=L^*/J$ can also be computed using $\aut L$. Let $S$ denote the
stabiliser in $\aut L$ of $J$ under the representation $\varrho$. Let $X$
denote a generating set for $S$. For each 
$\alpha\in X$ choose $\alpha^*\in\aut{L^*}$ as in the previous paragraph and let
$X^*=\{\alpha^*\ |\ \alpha\in X\}$. Suppose that $\{b_1,\ldots,b_d\}$ 
is a minimal generating set for $K$ and that $\{c_1,\ldots,c_l\}$ is a basis
for the last non-trivial term of the lower central series of $K$. 
For $i\in\{1,\ldots,d\}$ and $j\in\{1,\ldots,l\}$ let $\psi_{i,j}$ 
denote the automorphism that maps $b_i$ to $b_i+c_j$ 
and fixes $b_1,\ldots,b_{i-1},b_{i+1},\ldots,b_d$. Then 
$X^*\cup\{\psi_{i,j}\ |\ i=1,\ldots,d,\ 
j=1,\ldots,l\}$ is a generating set for $\aut K$. 

A similar approach to compute the automorphism group of a soluble Lie
algebra over a finite field is described in~\cite{eick}. Our method
is, however, more efficient for nilpotent Lie algebras.

The cover of a finite-dimensional 
nilpotent Lie algebra $L$ can be constructed in a way that is very 
similar  to the construction of the $p$-covering group 
of a finite $p$-group. A good description of this procedure can be found 
in~\cite{nnn}. 
Suppose that $L$ has class $c$, and hence the
lower central series is as follows:
$$
L= \lcs L1>\lcs L2=L'>\lcs L3>\cdots >\lcs Lc>\lcs L{c+1}=0.
$$
We say that  a basis $\mathcal B=\{b_1,\ldots,b_n\}$ for $L$ is 
{\em compatible with the
lower central series} if there are indices
$1=i_1<i_2<\cdots<i_{c-1}<i_c\leq n$ such that $\{b_{i_k},\ldots,b_n\}$ is
a basis of $\lcs Lk$ for $k\in\{1,\ldots,c\}$.

Suppose that $b_i\in\lcs L{j}\setminus\lcs L{j+1}$. Then we say that
the number $j$ is the {\em weight} of $b_i$. 
We call a basis $\mathcal B$ a {\em nilpotent basis} if the following
hold.
\begin{enumerate}
\item[(i)] The basis $\mathcal B$ is compatible with the lower central series;
\item[(ii)] for each $b_i\in \mathcal B$ with weight $w\geq 2$ there are
  $b_{j_1},\ b_{j_2}\in\mathcal B$ with weight 1 and $w-1$, respectively, such
  that $b_i=[b_{j_1},b_{j_2}]$. The product $[b_{j_1},b_{j_2}]$ is
  called the {\em definition} of $b_i$.
\end{enumerate}
If $\{b_1,\ldots,b_n\}$ is a nilpotent basis for a Lie algebra $L$, then there are coefficients $\alpha_{i,j}^k$ for $i<j<k$ such that 
\begin{equation}\label{sct}
[b_i,b_j]=\sum_{k=j+1}^n\alpha_{i,j}^k b_k.
\end{equation}

It is routine to see that every finitely generated 
nilpotent Lie algebra has a nilpotent basis.

Suppose that $\mathcal B=\{b_1,\ldots,b_n\}$ is a nilpotent basis for
a $d$-generator, nilpotent Lie algebra and the $\alpha_{i,j}^k$ are as in~\eqref{sct}. We build a presentation for the Lie algebra $L^*$ as follows. 
The set $\{b_{d+1},\ldots,b_n\}$ is a basis for $L'$. 
If, for some $i<j$,  
the product $[b_i,b_j]$ is not a definition and $w(b_i)+w(b_j)\leq c+1$, 
then we modify the product in~\eqref{sct}
by introducing a central basis element $b_{i,j}$ and set
$$
[b_i,b_j]=\sum_{k=j+1}^n\alpha_{i,j}^k b_k+b_{i,j}.
$$
We introduce the new basis elements so that different non-defining 
products $[b_i,b_j]$ are augmented with different basis elements $b_{i,j}$. 
We also ensure that the newly introduced basis elements $b_{i,j}$ are central.
If a product $[b_i,b_j]$ is a definition of $b_k$, say, then 
the product $[b_i,b_j]=b_k$ is not modified. Similarly if $w(b_i)+w(b_j)> c+1$
then $[b_i,b_j]$ is left untouched. 
This way we obtain an anti-commutative algebra $\hat L$ with basis 
$\{b_1,\ldots,b_d\}\cup\{b_{i,j}\}$ where the product of two basis elements is 
defined using the rules above.
We compute the ideal $J$ in $\hat L$ generated by the set of elements
$$
\{[b_i,b_j,b_k]+[b_j,b_k,b_i]+[b_k,b_i,b_j]\ |\ i,\ j,\ k\in\{1,\ldots,n\}\}.
$$
Then we obtain the cover $L^*$ as $\hat L/J$. 

It is possible to make this basic algorithm to compute the cover more 
effective. In practice we only introduce a new basis element for products of 
the form $[b_i,b_j]$ where $w(b_i)=1$ and compute products $[b_i,b_j]$ with
$w(b_i)>1$ using the Jacobi identity.
We also use the result  in~\cite{hnv-l}
that $J$ is already generated by the set of 
element
$$
\{[b_i,b_j,b_k]+[b_j,b_k,b_i]+[b_k,b_i,b_j]\ |\ i\in\{1,\ldots,d\},\ i<j<k\leq n\}.
$$

The proof that the resulting Lie algebra is isomorphic to $L^*$ is completely
analogous to that in the $p$-group case; see~\cite{nnn} for details.

\section{Some classifications of small Lie algebras}
\label{class}

In theory it is possible to use the procedures described in 
Section~\ref{coversec} 
to
classify nilpotent $\mathbb F_q$-Lie algebras of a given dimension
using a recursion. It is clear that there is a unique 1-dimensional
nilpotent Lie algebra over each field $\F_q$; the automorphism group of this
algebra is naturally isomorphic to the multiplicative group
$\mathbb F_q^*$.  Suppose that we have a complete
and irredundant list of nilpotent $\F_q$-Lie algebras of dimension
$1,\ldots,n-1$ for some $n\geq 2$ and we are also given the
automorphism groups of these algebras.
Up to isomorphism, there is exactly one $n$-dimensional abelian $\F_q$-Lie algebra.
Each non-abelian nilpotent Lie algebra with dimension $n$ 
is an immediate descendant of a smaller-dimensional Lie algebra. Hence, for each algebra $L$ with dimension $m$ in the precomputed  list we construct the Lie
cover $L^*$, the multiplicator $M$,  
and $\aut L\varrho$ where $\varrho$ is the representation in~\eqref{ro}. 
Then, using that $M$ is finite, we construct the
orbits of the $(m+\dim M-n)$-dimensional 
allowable subspaces under the finite linear group $\aut
L\varrho$. For each orbit representative $U$ we construct the quotient
$L^*/U$ and the stabiliser of $U$ in $\aut L$ under the representation
$\varrho$. The automorphism group of $L^*/U$ can now be constructed as
described in Section~\ref{coversec}. The collection of all Lie algebras
$L^*/U$ so obtained is a complete and irredundant list of the isomorphism 
types of the non-abelian 
nilpotent
Lie  algebras with dimension $n$. 

Suppose that $L$ is a finite-dimensional nilpotent Lie algebra. Let $c$ denote the class
of $L$. Then the {\em type} of the Lie algebra $L$ is the symbol
$$
[\dim L/L',\dim L'/\lcs L3,\ldots,\dim \lcs Lc][\dim Z(L)].
$$

It is well known that, over an arbitrary field, 
there is just one nilpotent Lie algebra with dimension
1 and 2. There are two nilpotent 
Lie algebras with dimension~3 (the types are $[3][3]$ 
and $[2,1][1]$), and 3~nilpotent Lie algebras with dimension~4 (the types are
$[4][4]$, $[3,1][2]$, $[2,1,1][1]$). The number of isomorphism types of 
5-dimensional nilpotent Lie algebras
is~9 over all fields; see~\cite{gozebook}. 
Up to isomorphism, there is exactly one Lie algebra with each of the following
types: $[5][5]$, $[4,1][3]$, $[4,1][1]$, $[3,2][2]$, $[3,1,1][2]$, $[3,1,1][1]$, $[2,1,2][2]$; there are two Lie algebras with type $[2,1,1,1][1]$. 

The number of 6-dimensional nilpotent Lie algebras depends 
on the underlying field. 
Using the {\sf GAP~4} package described in Section~\ref{implsec}, I obtained 36 
isomorphism classes of 6-dimensional nilpotent Lie algebras over $\F_2$, and
34~such classes over $\F_3$ and $\F_5$. It is mentioned in Wilkinson's 
paper~\cite{wilk} that the number of isomorphism classes of finite $p$-groups 
with order $p^6$ and exponent~$p$ is~34 whenever $p\geq 7$. 
Though there are several mistakes in the main
part of Wilkinson's paper, this particular 
claim has independently been checked by several people and
is widely considered to be true. 
Using the 
Lazard correspondence~\cite[Section~4]{obnv-l} we obtain that, for $p\geq 7$, there
are 34~pairwise non-isomorphic 6-dimensional 
nilpotent $\F_p$-Lie algebras. 
In fact, the computation referred to above implies that
this claim holds already for $p=3,\ 5$.
The number of 6-dimensional, nilpotent 
$\F_2$-Lie algebras for each possible type 
can be found in Table~1, 
while Table~2 contains the same information over $\F_3$ and $\F_5$. 
One can read off, for instance, from these tables that there are 6 pairwise
non-isomorphic nilpotent Lie algebras with type $[2,1,1,1,1][1]$ over $\F_2$
and there are only 5 such Lie algebras over $\F_3$ and $\F_5$.

It is reported in~\cite{gong} that~\cite{shed} contains a classification
of 6-dimensional nilpotent Lie algebras over any field, but this work is 
unpublished and contains several mistakes. 
There exist classifications of 6-dimensional nilpotent Lie algebras over
infinite fields; see for instance~\cite{gozebook}.

A classification of finite $p$-groups with exponent~7 and 
order $p^7$  was obtained
by Wilkinson~\cite{wilk}. If $p\geq 7$ then, by the Lazard correspondence,
the number of finite $p$-groups with exponent $p$ and order $p^7$ coincides 
with the number of 7-dimensional 
nilpotent $\F_p$-Lie algebras. 
According to Wilkinson this number is $173+7p+2\gcd(p-1,3)$, but Michael
Vaughan-Lee pointed out in private communication 
that there are several mistakes in Wilkinson's 
calculations and the correct number is
\begin{equation}\label{wilk}
174+7p+2\gcd(p-1,3).
\end{equation}
Computer calculations with the {\sf GAP~4} package 
described in Section~\ref{implsec} 
shows that the number of 7-dimensional 
nilpotent Lie algebras over $\F_2$, $\F_3$, and $\F_5$ is 202, 199, 211, 
respectively; the number of Lie algebras for each possible type is presented 
in Tables~3--5 in Section~\ref{tables}. This 
calculation also shows that~\eqref{wilk} is valid over $\F_5$. 
Michael Vaughan-Lee independently obtained a 
classification of nilpotent Lie rings with order $p^7$, and the numbers above
were also confirmed by his computation.

For some classifications of 7-dimensional nilpotent Lie algebras over infinite
fields we refer to~\cite{7C,7RC,gong,gozewww}

The author's {\sf GAP~4} program was also used the obtain a classification of 
nilpotent $\F_2$-Lie algebras with dimension 8 and 9. The total number of such 
Lie algebras is 1831 and 27073. More detailed information about the possible types
can be found in Tables~6--8 of Section~\ref{tables}.

The classifications of
nilpotent Lie algebras in Theorem~\ref{main} are 
available in {\sf GAP~4} format on 
the author's web site \verb+(www.sztaki.hu/~schneider/Research/SmallLie/)+.

\section{Implementation of the algorithms}\label{implsec}

Implementations
of all procedures described in Section~\ref{coversec} 
are available in the {\sf
 GAP~4} computer algebra package {\sf Sophus}. This program can freely
be downloaded from the author's web page \verb+(www.sztaki.hu/~schneider/Research/Sophus)+.
 The current version of {\sf Sophus} contains
\begin{itemize}
\item[(i)] a program to compute the cover of a nilpotent Lie algebra;
\item[(ii)] a program to compute the automorphism group of a nilpotent
  Lie
  algebra;
\item[(iii)] a program to compute the set of immediate descendants of
  a nilpotent Lie algebra;
\item[(iv)] a program to check if two nilpotent Lie algebras are
  isomorphic.
\end{itemize}
The full implementation of these procedures is nearly 4000 lines long.

The classifications presented in the previous section were 
computed on several Pentium~4 computers between 1.7 and 2.5 GHz CPU
speed and 1-2 GB memory. The computation of the list of the $\F_2$-Lie
algebras with dimension at most~6 takes only a few seconds, while those
of dimension~7 takes about 3 minutes.

Determining the remaining classes of nilpotent 
Lie algebras in Theorem~\ref{main} is more complicated and requires human
intervention. Most of the descendant computations for the 
8 and 9-dimensional Lie algebras over $\F_2$ could easily be 
carried out. However, computing the 
8-dimensional descendants of the 6-dimensional abelian Lie algebra 
requires finding representatives of 
the ${\sf GL}(6,2)$-orbits on the 
set of 178,940,587 allowable subspaces under
the action in~\eqref{ro}. In the computation of the 
9-dimensional descendants of the 7-dimensional abelian Lie algebra, the
number of allowable subspaces is 733,006,703,275.
In such cases I applied the Cauchy-Frobenius Lemma (see~\cite[Section~4]{o'be}) 
to predict the number of 
descendants.
Then I used either the ideas of O'Brien's extended algorithm 
presented in~\cite[Section~2]{256}, 
or the existing classification of $2$-groups of order at most 
$2^9$.
In the latter case
I constructed Lie algebras associated with the $2$-central series filtration 
of the groups, tested them for isomorphism, 
and eliminated the duplicates.

For 
computing the 7-dimensional 
descendants of the 5-dimensional abelian Lie algebras over $\F_3$ 
and $\F_5$, I used the result of the corresponding computation over $\F_2$. 
The Cauchy-Frobenius Lemma implies that the number of these Lie algebras is
the same over $\F_2$, $\F_3$, and $\F_5$. It is possible
to interpret the structure constants table of the $\F_2$-Lie algebras over
$\F_3$ and $\F_5$ and obtain the required lists. Then the algebras in these lists 
were tested for non-isomorphism.

The most difficult problem when computing the immediate 
descendants of a nilpotent Lie algebra is computing the orbits 
of the allowable subspaces under 
the representation~\eqref{ro}. Further, for computing the automorphism group of 
an immediate descendant, the stabiliser of an allowable subspace must also be
calculated; see Section~\ref{coversec}. These orbit-stabiliser computations were
carried out adopting the procedures described in~\cite{el-go'b}.

\section{The tables}\label{tables}

{\small
\begin{center}
\begin{equation*}\label{table6}
\begin{array}{cccc}\hline\begin{array}{c}{[ 6 ][ 6 ]} \\ 1\end{array}
&\begin{array}{c}{[ 5, 1 ][ 4 ]} \\ 1\end{array}
&\begin{array}{c}{[ 5, 1 ][ 2 ]} \\ 1\end{array}
&\begin{array}{c}{[ 4, 2 ][ 3 ]} \\ 1\end{array}
\\\hline
\begin{array}{c}{[ 4, 2 ][ 2 ]} \\ 3\end{array}
&\begin{array}{c}{[ 4, 1, 1 ][ 3 ]} \\ 1\end{array}
&\begin{array}{c}{[ 4, 1, 1 ][ 2 ]} \\ 1\end{array}
&\begin{array}{c}{[ 4, 1, 1 ][ 1 ]} \\ 1\end{array}
\\\hline
\begin{array}{c}{[ 3, 3 ][ 3 ]} \\ 1\end{array}
&\begin{array}{c}{[ 3, 2, 1 ][ 2 ]} \\ 3\end{array}
&\begin{array}{c}{[ 3, 2, 1 ][ 1 ]} \\ 3\end{array}
&\begin{array}{c}{[ 3, 1, 2 ][ 3 ]} \\ 1\end{array}
\\\hline
\begin{array}{c}{[ 3, 1, 2 ][ 2 ]} \\ 3\end{array}
&\begin{array}{c}{[ 3, 1, 1, 1 ][ 2 ]} \\ 2\end{array}
&\begin{array}{c}{[ 3, 1, 1, 1 ][ 1 ]} \\ 4\end{array}
&\begin{array}{c}{[ 2, 1, 2, 1 ][ 2 ]} \\ 1\end{array}
\\\hline
\begin{array}{c}{[ 2, 1, 2, 1 ][ 1 ]} \\ 2\end{array}
&\begin{array}{c}{[ [ 2, 1, 1, 1, 1 ], [ 1 ] ]} \\ 6\end{array}\\\hline\end{array}\end{equation*}
Table~1: The nilpotent Lie algebras of dimension 6 over $\F_2$
\end{center}

\begin{center}
$$
\begin{array}{cccc}\hline\begin{array}{c}{[ 6 ][ 6 ]} \\ 1\end{array}
&\begin{array}{c}{[ 5, 1 ][ 4 ]} \\ 1\end{array}
&\begin{array}{c}{[ 5, 1 ][ 2 ]} \\ 1\end{array}
&\begin{array}{c}{[ 4, 2 ][ 3 ]} \\ 1\end{array}
\\\hline
\begin{array}{c}{[ 4, 2 ][ 2 ]} \\ 3\end{array}
&\begin{array}{c}{[ 4, 1, 1 ][ 3 ]} \\ 1\end{array}
&\begin{array}{c}{[ 4, 1, 1 ][ 2 ]} \\ 1\end{array}
&\begin{array}{c}{[ 4, 1, 1 ][ 1 ]} \\ 1\end{array}
\\\hline
\begin{array}{c}{[ 3, 3 ][ 3 ]} \\ 1\end{array}
&\begin{array}{c}{[ 3, 2, 1 ][ 2 ]} \\ 3\end{array}
&\begin{array}{c}{[ 3, 2, 1 ][ 1 ]} \\ 3\end{array}
&\begin{array}{c}{[ 3, 1, 2 ][ 3 ]} \\ 1\end{array}
\\\hline
\begin{array}{c}{[ 3, 1, 2 ][ 2 ]} \\ 3\end{array}
&\begin{array}{c}{[ 3, 1, 1, 1 ][ 2 ]} \\ 2\end{array}
&\begin{array}{c}{[ 3, 1, 1, 1 ][ 1 ]} \\ 3\end{array}
&\begin{array}{c}{[ 2, 1, 2, 1 ][ 2 ]} \\ 1\end{array}
\\\hline
\begin{array}{c}{[ 2, 1, 2, 1 ][ 1 ]} \\ 2\end{array}
&\begin{array}{c}{[ [ 2, 1, 1, 1, 1 ], [ 1 ] ]} \\ 5\end{array}\\\hline\end{array}$$
Table~2: The nilpotent Lie algebras with dimension 6 over $\F_3$ and $\F_5$
\end{center}

\newpage

\begin{center}
$$
\begin{array}{cccc}\hline\begin{array}{c}{[ 7 ][ 7 ]} \\ 1\end{array}
&\begin{array}{c}{[ 6, 1 ][ 5 ]} \\ 1\end{array}
&\begin{array}{c}{[ 6, 1 ][ 3 ]} \\ 1\end{array}
&\begin{array}{c}{[ 6, 1 ][ 1 ]} \\ 1\end{array}
\\\hline
\begin{array}{c}{[ 5, 2 ][ 4 ]} \\ 1\end{array}
&\begin{array}{c}{[ 5, 2 ][ 3 ]} \\ 3\end{array}
&\begin{array}{c}{[ 5, 2 ][ 2 ]} \\ 2\end{array}
&\begin{array}{c}{[ 5, 1, 1 ][ 4 ]} \\ 1\end{array}
\\\hline
\begin{array}{c}{[ 5, 1, 1 ][ 3 ]} \\ 1\end{array}
&\begin{array}{c}{[ 5, 1, 1 ][ 2 ]} \\ 1\end{array}
&\begin{array}{c}{[ 5, 1, 1 ][ 1 ]} \\ 1\end{array}
&\begin{array}{c}{[ 4, 3 ][ 4 ]} \\ 1\end{array}
\\\hline
\begin{array}{c}{[ 4, 3 ][ 3 ]} \\ 5\end{array}
&\begin{array}{c}{[ 4, 2, 1 ][ 3 ]} \\ 3\end{array}
&\begin{array}{c}{[ 4, 2, 1 ][ 2 ]} \\ 12\end{array}
&\begin{array}{c}{[ 4, 2, 1 ][ 1 ]} \\ 9\end{array}
\\\hline
\begin{array}{c}{[ 4, 1, 2 ][ 4 ]} \\ 1\end{array}
&\begin{array}{c}{[ 4, 1, 2 ][ 3 ]} \\ 3\end{array}
&\begin{array}{c}{[ 4, 1, 2 ][ 2 ]} \\ 5\end{array}
&\begin{array}{c}{[ 4, 1, 1, 1 ][ 3 ]} \\ 2\end{array}
\\\hline
\begin{array}{c}{[ 4, 1, 1, 1 ][ 2 ]} \\ 4\end{array}
&\begin{array}{c}{[ 4, 1, 1, 1 ][ 1 ]} \\ 5\end{array}
&\begin{array}{c}{[ 3, 3, 1 ][ 3 ]} \\ 1\end{array}
&\begin{array}{c}{[ 3, 3, 1 ][ 2 ]} \\ 3\end{array}
\\\hline
\begin{array}{c}{[ 3, 3, 1 ][ 1 ]} \\ 2\end{array}
&\begin{array}{c}{[ 3, 2, 2 ][ 3 ]} \\ 2\end{array}
&\begin{array}{c}{[ 3, 2, 2 ][ 2 ]} \\ 21\end{array}
&\begin{array}{c}{[ 3, 2, 1, 1 ][ 2 ]} \\ 9\end{array}
\\\hline
\begin{array}{c}{[ 3, 2, 1, 1 ][ 1 ]} \\ 13\end{array}
&\begin{array}{c}{[ 3, 1, 2, 1 ][ 3 ]} \\ 1\end{array}
&\begin{array}{c}{[ 3, 1, 2, 1 ][ 2 ]} \\ 11\end{array}
&\begin{array}{c}{[ 3, 1, 2, 1 ][ 1 ]} \\ 8\end{array}
\\\hline
\begin{array}{c}{[ 3, 1, 1, 1, 1 ][ 2 ]} \\ 6\end{array}
&\begin{array}{c}{[ 3, 1, 1, 1, 1 ][ 1 ]} \\ 21\end{array}
&\begin{array}{c}{[ 2, 1, 2, 2 ][ 2 ]} \\ 3\end{array}
&\begin{array}{c}{[ 2, 1, 2, 1, 1 ][ 2 ]} \\ 4\end{array}
\\\hline
\begin{array}{c}{[ 2, 1, 2, 1, 1 ][ 1 ]} \\ 14\end{array}
&\begin{array}{c}{[ 2, 1, 1, 1, 2 ][ 2 ]} \\ 4\end{array}
&\begin{array}{c}{[ [ 2, 1, 1, 1, 1, 1 ], [ 1 ] ]} \\ 15\end{array}\\\hline\end{array}
$$
\label{table7}
Table~3: The nilpotent Lie algebras of dimension 7 over $\F_2$
\end{center}

\begin{center}
\vbox{$$
\begin{array}{cccc}\hline\begin{array}{c}{[ 7 ][ 7 ]} \\ 1\end{array}
&\begin{array}{c}{[ 6, 1 ][ 5 ]} \\ 1\end{array}
&\begin{array}{c}{[ 6, 1 ][ 3 ]} \\ 1\end{array}
&\begin{array}{c}{[ 6, 1 ][ 1 ]} \\ 1\end{array}
\\\hline
\begin{array}{c}{[ 5, 2 ][ 4 ]} \\ 1\end{array}
&\begin{array}{c}{[ 5, 2 ][ 3 ]} \\ 3\end{array}
&\begin{array}{c}{[ 5, 2 ][ 2 ]} \\ 2\end{array}
&\begin{array}{c}{[ 5, 1, 1 ][ 4 ]} \\ 1\end{array}
\\\hline
\begin{array}{c}{[ 5, 1, 1 ][ 3 ]} \\ 1\end{array}
&\begin{array}{c}{[ 5, 1, 1 ][ 2 ]} \\ 1\end{array}
&\begin{array}{c}{[ 5, 1, 1 ][ 1 ]} \\ 1\end{array}
&\begin{array}{c}{[ 4, 3 ][ 4 ]} \\ 1\end{array}
\\\hline
\begin{array}{c}{[ 4, 3 ][ 3 ]} \\ 5\end{array}
&\begin{array}{c}{[ 4, 2, 1 ][ 3 ]} \\ 3\end{array}
&\begin{array}{c}{[ 4, 2, 1 ][ 2 ]} \\ 12\end{array}
&\begin{array}{c}{[ 4, 2, 1 ][ 1 ]} \\ 9\end{array}
\\\hline
\begin{array}{c}{[ 4, 1, 2 ][ 4 ]} \\ 1\end{array}
&\begin{array}{c}{[ 4, 1, 2 ][ 3 ]} \\ 3\end{array}
&\begin{array}{c}{[ 4, 1, 2 ][ 2 ]} \\ 5\end{array}
&\begin{array}{c}{[ 4, 1, 1, 1 ][ 3 ]} \\ 2\end{array}
\\\hline
\begin{array}{c}{[ 4, 1, 1, 1 ][ 2 ]} \\ 3\end{array}
&\begin{array}{c}{[ 4, 1, 1, 1 ][ 1 ]} \\ 5\end{array}
&\begin{array}{c}{[ 3, 3, 1 ][ 3 ]} \\ 1\end{array}
&\begin{array}{c}{[ 3, 3, 1 ][ 2 ]} \\ 3\end{array}
\\\hline
\begin{array}{c}{[ 3, 3, 1 ][ 1 ]} \\ 5\end{array}
&\begin{array}{c}{[ 3, 2, 2 ][ 3 ]} \\ 2\end{array}
&\begin{array}{c}{[ 3, 2, 2 ][ 2 ]} \\ 21\end{array}
&\begin{array}{c}{[ 3, 2, 1, 1 ][ 2 ]} \\ 8\end{array}
\\\hline
\begin{array}{c}{[ 3, 2, 1, 1 ][ 1 ]} \\ 14\end{array}
&\begin{array}{c}{[ 3, 1, 2, 1 ][ 3 ]} \\ 1\end{array}
&\begin{array}{c}{[ 3, 1, 2, 1 ][ 2 ]} \\ 10\end{array}
&\begin{array}{c}{[ 3, 1, 2, 1 ][ 1 ]} \\ 12\end{array}
\\\hline
\begin{array}{c}{[ 3, 1, 1, 1, 1 ][ 2 ]} \\ 5\end{array}
&\begin{array}{c}{[ 3, 1, 1, 1, 1 ][ 1 ]} \\ 17\end{array}
&\begin{array}{c}{[ 2, 1, 2, 2 ][ 2 ]} \\ 3\end{array}
&\begin{array}{c}{[ 2, 1, 2, 1, 1 ][ 2 ]} \\ 4\end{array}
\\\hline
\begin{array}{c}{[ 2, 1, 2, 1, 1 ][ 1 ]} \\ 16\end{array}
&\begin{array}{c}{[ 2, 1, 1, 1, 2 ][ 2 ]} \\ 3\end{array}
&\begin{array}{c}{[ [ 2, 1, 1, 1, 1, 1 ], [ 1 ] ]} \\ 11\end{array}\\\hline\end{array}$$
Table~4: The nilpotent Lie algebras with dimension 7 over $\F_3$}
\end{center}
}

\begin{center}
\vbox{
$$
\begin{array}{cccc}\hline\begin{array}{c}{[ 7 ][ 7 ]} \\ 1\end{array}
&\begin{array}{c}{[ 6, 1 ][ 5 ]} \\ 1\end{array}
&\begin{array}{c}{[ 6, 1 ][ 3 ]} \\ 1\end{array}
&\begin{array}{c}{[ 6, 1 ][ 1 ]} \\ 1\end{array}
\\\hline
\begin{array}{c}{[ 5, 2 ][ 4 ]} \\ 1\end{array}
&\begin{array}{c}{[ 5, 2 ][ 3 ]} \\ 3\end{array}
&\begin{array}{c}{[ 5, 2 ][ 2 ]} \\ 2\end{array}
&\begin{array}{c}{[ 5, 1, 1 ][ 4 ]} \\ 1\end{array}
\\\hline
\begin{array}{c}{[ 5, 1, 1 ][ 3 ]} \\ 1\end{array}
&\begin{array}{c}{[ 5, 1, 1 ][ 2 ]} \\ 1\end{array}
&\begin{array}{c}{[ 5, 1, 1 ][ 1 ]} \\ 1\end{array}
&\begin{array}{c}{[ 4, 3 ][ 4 ]} \\ 1\end{array}
\\\hline
\begin{array}{c}{[ 4, 3 ][ 3 ]} \\ 5\end{array}
&\begin{array}{c}{[ 4, 2, 1 ][ 3 ]} \\ 3\end{array}
&\begin{array}{c}{[ 4, 2, 1 ][ 2 ]} \\ 12\end{array}
&\begin{array}{c}{[ 4, 2, 1 ][ 1 ]} \\ 9\end{array}
\\\hline
\begin{array}{c}{[ 4, 1, 2 ][ 4 ]} \\ 1\end{array}
&\begin{array}{c}{[ 4, 1, 2 ][ 3 ]} \\ 3\end{array}
&\begin{array}{c}{[ 4, 1, 2 ][ 2 ]} \\ 5\end{array}
&\begin{array}{c}{[ 4, 1, 1, 1 ][ 3 ]} \\ 2\end{array}
\\\hline
\begin{array}{c}{[ 4, 1, 1, 1 ][ 2 ]} \\ 3\end{array}
&\begin{array}{c}{[ 4, 1, 1, 1 ][ 1 ]} \\ 5\end{array}
&\begin{array}{c}{[ 3, 3, 1 ][ 3 ]} \\ 1\end{array}
&\begin{array}{c}{[ 3, 3, 1 ][ 2 ]} \\ 3\end{array}
\\\hline
\begin{array}{c}{[ 3, 3, 1 ][ 1 ]} \\ 6\end{array}
&\begin{array}{c}{[ 3, 2, 2 ][ 3 ]} \\ 2\end{array}
&\begin{array}{c}{[ 3, 2, 2 ][ 2 ]} \\ 21\end{array}
&\begin{array}{c}{[ 3, 2, 1, 1 ][ 2 ]} \\ 8\end{array}
\\\hline
\begin{array}{c}{[ 3, 2, 1, 1 ][ 1 ]} \\ 18\end{array}
&\begin{array}{c}{[ 3, 1, 2, 1 ][ 3 ]} \\ 1\end{array}
&\begin{array}{c}{[ 3, 1, 2, 1 ][ 2 ]} \\ 10\end{array}
&\begin{array}{c}{[ 3, 1, 2, 1 ][ 1 ]} \\ 16\end{array}
\\\hline
\begin{array}{c}{[ 3, 1, 1, 1, 1 ][ 2 ]} \\ 5\end{array}
&\begin{array}{c}{[ 3, 1, 1, 1, 1 ][ 1 ]} \\ 16\end{array}
&\begin{array}{c}{[ 2, 1, 2, 2 ][ 2 ]} \\ 3\end{array}
&\begin{array}{c}{[ 2, 1, 2, 1, 1 ][ 2 ]} \\ 4\end{array}
\\\hline
\begin{array}{c}{[ 2, 1, 2, 1, 1 ][ 1 ]} \\ 18\end{array}
&\begin{array}{c}{[ 2, 1, 1, 1, 2 ][ 2 ]} \\ 3\end{array}
&\begin{array}{c}{[ [ 2, 1, 1, 1, 1, 1 ], [ 1 ] ]} \\ 13\end{array}\\\hline\end{array}$$
Table 5: The nilpotent Lie algebras with dimension 7 over $\F_5$}
\end{center}

{\small
\begin{center}
$$
\begin{array}{cccc}\hline\begin{array}{c}{[ 8 ][ 8 ]} \\ 1\end{array}
&\begin{array}{c}{[ 7, 1 ][ 6 ]} \\ 1\end{array}
&\begin{array}{c}{[ 7, 1 ][ 4 ]} \\ 1\end{array}
&\begin{array}{c}{[ 7, 1 ][ 2 ]} \\ 1\end{array}
\\\hline
\begin{array}{c}{[ 6, 2 ][ 5 ]} \\ 1\end{array}
&\begin{array}{c}{[ 6, 2 ][ 4 ]} \\ 3\end{array}
&\begin{array}{c}{[ 6, 2 ][ 3 ]} \\ 2\end{array}
&\begin{array}{c}{[ 6, 2 ][ 2 ]} \\ 8\end{array}
\\\hline
\begin{array}{c}{[ 6, 1, 1 ][ 5 ]} \\ 1\end{array}
&\begin{array}{c}{[ 6, 1, 1 ][ 4 ]} \\ 1\end{array}
&\begin{array}{c}{[ 6, 1, 1 ][ 3 ]} \\ 1\end{array}
&\begin{array}{c}{[ 6, 1, 1 ][ 2 ]} \\ 1\end{array}
\\\hline
\begin{array}{c}{[ 6, 1, 1 ][ 1 ]} \\ 1\end{array}
&\begin{array}{c}{[ 5, 3 ][ 5 ]} \\ 1\end{array}
&\begin{array}{c}{[ 5, 3 ][ 4 ]} \\ 5\end{array}
&\begin{array}{c}{[ 5, 3 ][ 3 ]} \\ 16\end{array}
\\\hline
\begin{array}{c}{[ 5, 2, 1 ][ 4 ]} \\ 3\end{array}
&\begin{array}{c}{[ 5, 2, 1 ][ 3 ]} \\ 12\end{array}
&\begin{array}{c}{[ 5, 2, 1 ][ 2 ]} \\ 35\end{array}
&\begin{array}{c}{[ 5, 2, 1 ][ 1 ]} \\ 13\end{array}
\\\hline
\begin{array}{c}{[ 5, 1, 2 ][ 5 ]} \\ 1\end{array}
&\begin{array}{c}{[ 5, 1, 2 ][ 4 ]} \\ 3\end{array}
&\begin{array}{c}{[ 5, 1, 2 ][ 3 ]} \\ 5\end{array}
&\begin{array}{c}{[ 5, 1, 2 ][ 2 ]} \\ 14\end{array}
\\\hline
\begin{array}{c}{[ 5, 1, 1, 1 ][ 4 ]} \\ 2\end{array}
&\begin{array}{c}{[ 5, 1, 1, 1 ][ 3 ]} \\ 4\end{array}
&\begin{array}{c}{[ 5, 1, 1, 1 ][ 2 ]} \\ 5\end{array}
&\begin{array}{c}{[ 5, 1, 1, 1 ][ 1 ]} \\ 5\end{array}
\\\hline
\begin{array}{c}{[ 4, 4 ][ 4 ]} \\ 4\end{array}
&\begin{array}{c}{[ 4, 3, 1 ][ 4 ]} \\ 1\end{array}
&\begin{array}{c}{[ 4, 3, 1 ][ 3 ]} \\ 29\end{array}
&\begin{array}{c}{[ 4, 3, 1 ][ 2 ]} \\ 51\end{array}
\\\hline
\begin{array}{c}{[ 4, 3, 1 ][ 1 ]} \\ 25\end{array}
&\begin{array}{c}{[ 4, 2, 2 ][ 4 ]} \\ 2\end{array}
&\begin{array}{c}{[ 4, 2, 2 ][ 3 ]} \\ 48\end{array}
&\begin{array}{c}{[ 4, 2, 2 ][ 2 ]} \\ 209\end{array}
\\\hline
\begin{array}{c}{[ 4, 2, 1, 1 ][ 3 ]} \\ 9\end{array}
&\begin{array}{c}{[ 4, 2, 1, 1 ][ 2 ]} \\ 59\end{array}
&\begin{array}{c}{[ 4, 2, 1, 1 ][ 1 ]} \\ 54\end{array}
&\begin{array}{c}{[ 4, 1, 2, 1 ][ 4 ]} \\ 1\end{array}
\\\hline
\begin{array}{c}{[ 4, 1, 2, 1 ][ 3 ]} \\ 11\end{array}
&\begin{array}{c}{[ 4, 1, 2, 1 ][ 2 ]} \\ 48\end{array}
&\begin{array}{c}{[ 4, 1, 2, 1 ][ 1 ]} \\ 26\end{array}
&\begin{array}{c}{[ 4, 1, 1, 1, 1 ][ 3 ]} \\ 6\end{array}
\\\hline
\begin{array}{c}{[ 4, 1, 1, 1, 1 ][ 2 ]} \\ 21\end{array}
&\begin{array}{c}{[ 4, 1, 1, 1, 1 ][ 1 ]} \\ 39\end{array}
&\begin{array}{c}{[ 3, 3, 2 ][ 4 ]} \\ 1\end{array}
&\begin{array}{c}{[ 3, 3, 2 ][ 3 ]} \\ 15\end{array}
\\\hline
\begin{array}{c}{[ 3, 3, 2 ][ 2 ]} \\ 77\end{array}
&\begin{array}{c}{[ 3, 3, 1, 1 ][ 3 ]} \\ 3\end{array}
&\begin{array}{c}{[ 3, 3, 1, 1 ][ 2 ]} \\ 13\end{array}
&\begin{array}{c}{[ 3, 3, 1, 1 ][ 1 ]} \\ 6\end{array}
\\\hline
\begin{array}{c}{[ 3, 2, 3 ][ 3 ]} \\ 28\end{array}
&\begin{array}{c}{[ 3, 2, 2, 1 ][ 3 ]} \\ 11\end{array}
&\begin{array}{c}{[ 3, 2, 2, 1 ][ 2 ]} \\ 164\end{array}
&\begin{array}{c}{[ 3, 2, 2, 1 ][ 1 ]} \\ 84\end{array}
\\\hline
\begin{array}{c}{[ 3, 2, 1, 1, 1 ][ 2 ]} \\ 49\end{array}
&\begin{array}{c}{[ 3, 2, 1, 1, 1 ][ 1 ]} \\ 88\end{array}
&\begin{array}{c}{[ 3, 1, 2, 2 ][ 3 ]} \\ 3\end{array}
&\begin{array}{c}{[ 3, 1, 2, 2 ][ 2 ]} \\ 37\end{array}
\\\hline
\begin{array}{c}{[ 3, 1, 2, 1, 1 ][ 3 ]} \\ 4\end{array}
&\begin{array}{c}{[ 3, 1, 2, 1, 1 ][ 2 ]} \\ 71\end{array}
&\begin{array}{c}{[ 3, 1, 2, 1, 1 ][ 1 ]} \\ 82\end{array}
&\begin{array}{c}{[ 3, 1, 1, 1, 2 ][ 3 ]} \\ 4\end{array}
\\\hline
\begin{array}{c}{[ 3, 1, 1, 1, 2 ][ 2 ]} \\ 39\end{array}
&\begin{array}{c}{[ 3, 1, 1, 1, 1, 1 ][ 2 ]} \\ 15\end{array}
&\begin{array}{c}{[ 3, 1, 1, 1, 1, 1 ][ 1 ]} \\ 80\end{array}
&\begin{array}{c}{[ 2, 1, 2, 3 ][ 3 ]} \\ 1\end{array}
\\\hline
\begin{array}{c}{[ 2, 1, 2, 2, 1 ][ 2 ]} \\ 26\end{array}
&\begin{array}{c}{[ 2, 1, 2, 2, 1 ][ 1 ]} \\ 20\end{array}
&\begin{array}{c}{[ 2, 1, 2, 1, 2 ][ 3 ]} \\ 2\end{array}
&\begin{array}{c}{[ 2, 1, 2, 1, 2 ][ 2 ]} \\ 24\end{array}
\\\hline
\begin{array}{c}{[ 2, 1, 2, 1, 1, 1 ][ 2 ]} \\ 12\end{array}
&\begin{array}{c}{[ 2, 1, 2, 1, 1, 1 ][ 1 ]} \\ 24\end{array}
&\begin{array}{c}{[ 2, 1, 1, 1, 2, 1 ][ 2 ]} \\ 11\end{array}
&\begin{array}{c}{[ [ 2, 1, 1, 1, 1, 1, 1 ], [ 1 ] ]} \\ 47\end{array}\\\hline\end{array}$$
\label{table8}
Table 6: The nilpotent Lie algebras with dimension 8 over $\F_2$
\end{center}
}

{\small\begin{center}
$$
\begin{array}{cccc}\hline\begin{array}{c}{[ 9 ][ 9 ]} \\ 1\end{array}
&\begin{array}{c}{[ 8, 1 ][ 7 ]} \\ 1\end{array}
&\begin{array}{c}{[ 8, 1 ][ 5 ]} \\ 1\end{array}
&\begin{array}{c}{[ 8, 1 ][ 3 ]} \\ 1\end{array}
\\\hline
\begin{array}{c}{[ 8, 1 ][ 1 ]} \\ 1\end{array}
&\begin{array}{c}{[ 7, 2 ][ 6 ]} \\ 1\end{array}
&\begin{array}{c}{[ 7, 2 ][ 5 ]} \\ 3\end{array}
&\begin{array}{c}{[ 7, 2 ][ 4 ]} \\ 2\end{array}
\\\hline
\begin{array}{c}{[ 7, 2 ][ 3 ]} \\ 8\end{array}
&\begin{array}{c}{[ 7, 2 ][ 2 ]} \\ 6\end{array}
&\begin{array}{c}{[ 7, 1, 1 ][ 6 ]} \\ 1\end{array}
&\begin{array}{c}{[ 7, 1, 1 ][ 5 ]} \\ 1\end{array}
\\\hline
\begin{array}{c}{[ 7, 1, 1 ][ 4 ]} \\ 1\end{array}
&\begin{array}{c}{[ 7, 1, 1 ][ 3 ]} \\ 1\end{array}
&\begin{array}{c}{[ 7, 1, 1 ][ 2 ]} \\ 1\end{array}
&\begin{array}{c}{[ 7, 1, 1 ][ 1 ]} \\ 1\end{array}
\\\hline
\begin{array}{c}{[ 6, 3 ][ 6 ]} \\ 1\end{array}
&\begin{array}{c}{[ 6, 3 ][ 5 ]} \\ 5\end{array}
&\begin{array}{c}{[ 6, 3 ][ 4 ]} \\ 16\end{array}
&\begin{array}{c}{[ 6, 3 ][ 3 ]} \\ 122\end{array}
\\\hline
\begin{array}{c}{[ 6, 2, 1 ][ 5 ]} \\ 3\end{array}
&\begin{array}{c}{[ 6, 2, 1 ][ 4 ]} \\ 12\end{array}
&\begin{array}{c}{[ 6, 2, 1 ][ 3 ]} \\ 35\end{array}
&\begin{array}{c}{[ 6, 2, 1 ][ 2 ]} \\ 70\end{array}
\\\hline
\begin{array}{c}{[ 6, 2, 1 ][ 1 ]} \\ 18\end{array}
&\begin{array}{c}{[ 6, 1, 2 ][ 6 ]} \\ 1\end{array}
&\begin{array}{c}{[ 6, 1, 2 ][ 5 ]} \\ 3\end{array}
&\begin{array}{c}{[ 6, 1, 2 ][ 4 ]} \\ 5\end{array}
\\\hline
\begin{array}{c}{[ 6, 1, 2 ][ 3 ]} \\ 14\end{array}
&\begin{array}{c}{[ 6, 1, 2 ][ 2 ]} \\ 25\end{array}
&\begin{array}{c}{[ 6, 1, 1, 1 ][ 5 ]} \\ 2\end{array}
&\begin{array}{c}{[ 6, 1, 1, 1 ][ 4 ]} \\ 4\end{array}
\\\hline
\begin{array}{c}{[ 6, 1, 1, 1 ][ 3 ]} \\ 5\end{array}
&\begin{array}{c}{[ 6, 1, 1, 1 ][ 2 ]} \\ 5\end{array}
&\begin{array}{c}{[ 6, 1, 1, 1 ][ 1 ]} \\ 5\end{array}
&\begin{array}{c}{[ 5, 4 ][ 5 ]} \\ 4\end{array}
\\\hline
\begin{array}{c}{[ 5, 4 ][ 4 ]} \\ 53\end{array}
&\begin{array}{c}{[ 5, 3, 1 ][ 5 ]} \\ 1\end{array}
&\begin{array}{c}{[ 5, 3, 1 ][ 4 ]} \\ 29\end{array}
&\begin{array}{c}{[ 5, 3, 1 ][ 3 ]} \\ 327\end{array}
\\\hline
\begin{array}{c}{[ 5, 3, 1 ][ 2 ]} \\ 318\end{array}
&\begin{array}{c}{[ 5, 3, 1 ][ 1 ]} \\ 133\end{array}
&\begin{array}{c}{[ 5, 2, 2 ][ 5 ]} \\ 2\end{array}
&\begin{array}{c}{[ 5, 2, 2 ][ 4 ]} \\ 48\end{array}
\\\hline
\begin{array}{c}{[ 5, 2, 2 ][ 3 ]} \\ 502\end{array}
&\begin{array}{c}{[ 5, 2, 2 ][ 2 ]} \\ 799\end{array}
&\begin{array}{c}{[ 5, 2, 1, 1 ][ 4 ]} \\ 9\end{array}
&\begin{array}{c}{[ 5, 2, 1, 1 ][ 3 ]} \\ 59\end{array}
\\\hline
\begin{array}{c}{[ 5, 2, 1, 1 ][ 2 ]} \\ 231\end{array}
&\begin{array}{c}{[ 5, 2, 1, 1 ][ 1 ]} \\ 129\end{array}
&\begin{array}{c}{[ 5, 1, 2, 1 ][ 5 ]} \\ 1\end{array}
&\begin{array}{c}{[ 5, 1, 2, 1 ][ 4 ]} \\ 11\end{array}
\\\hline
\begin{array}{c}{[ 5, 1, 2, 1 ][ 3 ]} \\ 48\end{array}
&\begin{array}{c}{[ 5, 1, 2, 1 ][ 2 ]} \\ 180\end{array}
&\begin{array}{c}{[ 5, 1, 2, 1 ][ 1 ]} \\ 37\end{array}
&\begin{array}{c}{[ 5, 1, 1, 1, 1 ][ 4 ]} \\ 6\end{array}
\\\hline
\begin{array}{c}{[ 5, 1, 1, 1, 1 ][ 3 ]} \\ 21\end{array}
&\begin{array}{c}{[ 5, 1, 1, 1, 1 ][ 2 ]} \\ 39\end{array}
&\begin{array}{c}{[ 5, 1, 1, 1, 1 ][ 1 ]} \\ 47\end{array}
&\begin{array}{c}{[ 4, 5 ][ 5 ]} \\ 2\end{array}
\\\hline
\begin{array}{c}{[ 4, 4, 1 ][ 4 ]} \\ 19\end{array}
&\begin{array}{c}{[ 4, 4, 1 ][ 3 ]} \\ 77\end{array}
&\begin{array}{c}{[ 4, 4, 1 ][ 2 ]} \\ 127\end{array}
&\begin{array}{c}{[ 4, 4, 1 ][ 1 ]} \\ 54\end{array}
\\\hline
\begin{array}{c}{[ 4, 3, 2 ][ 5 ]} \\ 1\end{array}
&\begin{array}{c}{[ 4, 3, 2 ][ 4 ]} \\ 55\end{array}
&\begin{array}{c}{[ 4, 3, 2 ][ 3 ]} \\ 814\end{array}
&\begin{array}{c}{[ 4, 3, 2 ][ 2 ]} \\ 2510\end{array}
\\\hline
\begin{array}{c}{[ 4, 3, 1, 1 ][ 4 ]} \\ 3\end{array}
&\begin{array}{c}{[ 4, 3, 1, 1 ][ 3 ]} \\ 131\end{array}
&\begin{array}{c}{[ 4, 3, 1, 1 ][ 2 ]} \\ 396\end{array}
&\begin{array}{c}{[ 4, 3, 1, 1 ][ 1 ]} \\ 296\end{array}
\\\hline
\begin{array}{c}{[ 4, 2, 3 ][ 4 ]} \\ 28\end{array}
&\begin{array}{c}{[ 4, 2, 3 ][ 3 ]} \\ 1377\end{array}
&\begin{array}{c}{[ 4, 2, 2, 1 ][ 4 ]} \\ 11\end{array}
&\begin{array}{c}{[ 4, 2, 2, 1 ][ 3 ]} \\ 402\end{array}
\\\hline
\end{array}
$$
Table~7: The nilpotent Lie algebras with dimension 9 over $\F_2$ (continued on
the next page)

$$
\begin{array}{cccc}
\hline
\begin{array}{c}{[ 4, 2, 2, 1 ][ 2 ]} \\ 2859\end{array}
&\begin{array}{c}{[ 4, 2, 2, 1 ][ 1 ]} \\ 713\end{array}
&\begin{array}{c}{[ 4, 2, 1, 1, 1 ][ 3 ]} \\ 49\end{array}
&\begin{array}{c}{[ 4, 2, 1, 1, 1 ][ 2 ]} \\ 487\end{array}
\\\hline
\begin{array}{c}{[ 4, 2, 1, 1, 1 ][ 1 ]} \\ 565\end{array}
&\begin{array}{c}{[ 4, 1, 2, 2 ][ 4 ]} \\ 3\end{array}
&\begin{array}{c}{[ 4, 1, 2, 2 ][ 3 ]} \\ 37\end{array}
&\begin{array}{c}{[ 4, 1, 2, 2 ][ 2 ]} \\ 258\end{array}
\\\hline
\begin{array}{c}{[ 4, 1, 2, 1, 1 ][ 4 ]} \\ 4\end{array}
&\begin{array}{c}{[ 4, 1, 2, 1, 1 ][ 3 ]} \\ 71\end{array}
&\begin{array}{c}{[ 4, 1, 2, 1, 1 ][ 2 ]} \\ 463\end{array}
&\begin{array}{c}{[ 4, 1, 2, 1, 1 ][ 1 ]} \\ 318\end{array}
\\\hline
\begin{array}{c}{[ 4, 1, 1, 1, 2 ][ 4 ]} \\ 4\end{array}
&\begin{array}{c}{[ 4, 1, 1, 1, 2 ][ 3 ]} \\ 39\end{array}
&\begin{array}{c}{[ 4, 1, 1, 1, 2 ][ 2 ]} \\ 191\end{array}
&\begin{array}{c}{[ 4, 1, 1, 1, 1, 1 ][ 3 ]} \\ 15\end{array}
\\\hline
\begin{array}{c}{[ 4, 1, 1, 1, 1, 1 ][ 2 ]} \\ 80\end{array}
&\begin{array}{c}{[ 4, 1, 1, 1, 1, 1 ][ 1 ]} \\ 213\end{array}
&\begin{array}{c}{[ 3, 3, 3 ][ 4 ]} \\ 16\end{array}
&\begin{array}{c}{[ 3, 3, 3 ][ 3 ]} \\ 642\end{array}
\\\hline
\begin{array}{c}{[ 3, 3, 2, 1 ][ 4 ]} \\ 2\end{array}
&\begin{array}{c}{[ 3, 3, 2, 1 ][ 3 ]} \\ 104\end{array}
&\begin{array}{c}{[ 3, 3, 2, 1 ][ 2 ]} \\ 808\end{array}
&\begin{array}{c}{[ 3, 3, 2, 1 ][ 1 ]} \\ 316\end{array}
\\\hline
\begin{array}{c}{[ 3, 3, 1, 1, 1 ][ 3 ]} \\ 16\end{array}
&\begin{array}{c}{[ 3, 3, 1, 1, 1 ][ 2 ]} \\ 86\end{array}
&\begin{array}{c}{[ 3, 3, 1, 1, 1 ][ 1 ]} \\ 76\end{array}
&\begin{array}{c}{[ 3, 2, 4 ][ 4 ]} \\ 12\end{array}
\\\hline
\begin{array}{c}{[ 3, 2, 3, 1 ][ 3 ]} \\ 258\end{array}
&\begin{array}{c}{[ 3, 2, 3, 1 ][ 2 ]} \\ 429\end{array}
&\begin{array}{c}{[ 3, 2, 3, 1 ][ 1 ]} \\ 203\end{array}
&\begin{array}{c}{[ 3, 2, 2, 2 ][ 3 ]} \\ 44\end{array}
\\\hline
\begin{array}{c}{[ 3, 2, 2, 2 ][ 2 ]} \\ 908\end{array}
&\begin{array}{c}{[ 3, 2, 2, 1, 1 ][ 3 ]} \\ 71\end{array}
&\begin{array}{c}{[ 3, 2, 2, 1, 1 ][ 2 ]} \\ 1296\end{array}
&\begin{array}{c}{[ 3, 2, 2, 1, 1 ][ 1 ]} \\ 1282\end{array}
\\\hline
\begin{array}{c}{[ 3, 2, 1, 1, 2 ][ 3 ]} \\ 33\end{array}
&\begin{array}{c}{[ 3, 2, 1, 1, 2 ][ 2 ]} \\ 325\end{array}
&\begin{array}{c}{[ 3, 2, 1, 1, 1, 1 ][ 2 ]} \\ 163\end{array}
&\begin{array}{c}{[ 3, 2, 1, 1, 1, 1 ][ 1 ]} \\ 435\end{array}
\\\hline
\begin{array}{c}{[ 3, 1, 2, 3 ][ 4 ]} \\ 1\end{array}
&\begin{array}{c}{[ 3, 1, 2, 3 ][ 3 ]} \\ 21\end{array}
&\begin{array}{c}{[ 3, 1, 2, 2, 1 ][ 3 ]} \\ 26\end{array}
&\begin{array}{c}{[ 3, 1, 2, 2, 1 ][ 2 ]} \\ 622\end{array}
\\\hline
\begin{array}{c}{[ 3, 1, 2, 2, 1 ][ 1 ]} \\ 302\end{array}
&\begin{array}{c}{[ 3, 1, 2, 1, 2 ][ 4 ]} \\ 2\end{array}
&\begin{array}{c}{[ 3, 1, 2, 1, 2 ][ 3 ]} \\ 79\end{array}
&\begin{array}{c}{[ 3, 1, 2, 1, 2 ][ 2 ]} \\ 353\end{array}
\\\hline
\begin{array}{c}{[ 3, 1, 2, 1, 1, 1 ][ 3 ]} \\ 12\end{array}
&\begin{array}{c}{[ 3, 1, 2, 1, 1, 1 ][ 2 ]} \\ 230\end{array}
&\begin{array}{c}{[ 3, 1, 2, 1, 1, 1 ][ 1 ]} \\ 314\end{array}
&\begin{array}{c}{[ 3, 1, 1, 1, 2, 1 ][ 3 ]} \\ 11\end{array}
\\\hline
\begin{array}{c}{[ 3, 1, 1, 1, 2, 1 ][ 2 ]} \\ 181\end{array}
&\begin{array}{c}{[ 3, 1, 1, 1, 1, 1, 1 ][ 2 ]} \\ 47\end{array}
&\begin{array}{c}{[ 3, 1, 1, 1, 1, 1, 1 ][ 1 ]} \\ 423\end{array}
&\begin{array}{c}{[ 2, 1, 2, 3, 1 ][ 3 ]} \\ 5\end{array}
\\\hline
\begin{array}{c}{[ 2, 1, 2, 3, 1 ][ 2 ]} \\ 10\end{array}
&\begin{array}{c}{[ 2, 1, 2, 2, 2 ][ 3 ]} \\ 19\end{array}
&\begin{array}{c}{[ 2, 1, 2, 2, 2 ][ 2 ]} \\ 170\end{array}
&\begin{array}{c}{[ 2, 1, 2, 2, 1, 1 ][ 2 ]} \\ 60\end{array}
\\\hline
\begin{array}{c}{[ 2, 1, 2, 2, 1, 1 ][ 1 ]} \\ 98\end{array}
&\begin{array}{c}{[ 2, 1, 2, 1, 2, 1 ][ 3 ]} \\ 6\end{array}
&\begin{array}{c}{[ 2, 1, 2, 1, 2, 1 ][ 2 ]} \\ 62\end{array}
&\begin{array}{c}{[ 2, 1, 2, 1, 2, 1 ][ 1 ]} \\ 16\end{array}
\\\hline
\begin{array}{c}{[ 2, 1, 2, 1, 1, 1, 1 ][ 2 ]} \\ 40\end{array}
&\begin{array}{c}{[ 2, 1, 2, 1, 1, 1, 1 ][ 1 ]} \\ 124\end{array}
&\begin{array}{c}{[ 2, 1, 1, 1, 2, 2 ][ 2 ]} \\ 7\end{array}
&\begin{array}{c}{[ 2, 1, 1, 1, 2, 1, 1 ][ 2 ]} \\ 45\end{array}
\\\hline
\begin{array}{c}{[ 2, 1, 1, 1, 2, 1, 1 ][ 1 ]} \\ 18\end{array}
&\begin{array}{c}{[ 2, 1, 1, 1, 1, 1, 2 ][ 2 ]} \\ 32\end{array}
&\begin{array}{c}{[ [ 2, 1, 1, 1, 1, 1, 1, 1 ], [ 1 ] ]} \\ 
124\end{array}\\\hline\end{array}$$
\label{table9}
Table~8: The nilpotent Lie algebras with dimension 9 over $\F_2$ 
(continued from the previous page)
\end{center}
}

\section{Acknowledgments}

A large part of the research presented in this paper was carried out at 
the Technische Universit\"at Braunschweig. I am grateful for Bettina Eick for
her interest in this project. I am also indebted to Michael Vaughan-Lee for 
sharing his work concerning the 7-dimensional Lie algebras; to Marco 
Costantini for testing my results
and finding a bug in my descendant computation;
and to J\"org Feldvoss for reading an earlier draft.

\end{document}